# Degree Conditions for Dominating Cycles in 1-tough Graphs

Zh.G. Nikoghosyan*

February 13, 2012

**Abstract**

We prove: (i) if $G$ is a 1-tough graph of order $n$ and minimum degree $\delta$ with $\delta \geq (n-2)/3$ then each longest cycle in $G$ is a dominating cycle unless $G$ belongs to an easily specified class of graphs with $\kappa(G) = 2$ and $\tau(G) = 1$. The second result follows immediately from the first result: (ii) if $G$ is a 3-connected 1-tough graph with $\delta \geq (n-2)/3$ then each longest cycle in $G$ is a dominating cycle.

## 1 Introduction

Only finite undirected graphs without loops or multiple edges are considered. We reserve $n$, $\delta$, $\kappa$ and $\tau$ to denote the number of vertices (order), the minimum degree, connectivity and the toughness of a graph, respectively. A good reference for any undefined terms is [3].

The earliest sufficient condition for dominating cycles was developed in 1971 due to Nash-Williams [4].

**Theorem A [4].** Let $G$ be a 2-connected graph. If $\delta \geq \frac{1}{3}(n+2)$ then each longest cycle in $G$ is a dominating cycle.

In 1979, Bigalke and Jung [2] proved that the minimum degree condition $\delta \geq (n+2)/3$ in Theorem A can be slightly relaxed by replacing the 2-connectivity condition with stronger 1-tough condition.

**Theorem B [2].** Let $G$ be a 1-tough graph. If $\delta \geq \frac{1}{3}n$ then each longest cycle in $G$ is a dominating cycle.

---

*G.G. Nicoghossian (up to 1997)

Recently it was proved that the minimum degree bound $n/3$ in Theorem B can be lowered to $(n-2)/3$ when $\tau > 1$.

**Theorem C** [5]. Let $G$ be a graph with $\tau > 1$. If $\delta \geq \frac{1}{3}(n-2)$ then each longest cycle in $G$ is a dominating cycle.

In this paper we prove two analogous theorems concerning 1-tough graphs with $\delta \geq \frac{1}{3}(n-2)$. For this purpose, we need to describe some classes of graphs. For two given graphs $H_1$ and $H_2$, we write $H_1 \subseteq H_2$ if $V(H_1) \subseteq V(H_2)$ and $E(H_1) \subseteq E(H_2)$.

Let $H_1, H_2, H_3$ be disjoint connected graphs with $x_i, y_i \in V(H_i)$ ($i = 1, 2, 3$) and $d \geq 3$ an integer. Form a graph $H$ by identifying $x_1, x_2, x_3$ and adding the edges of a triangle between $y_1, y_2, y_3$. We will say that $H \in \Re_1$ if and only if

$$|V(H_1)| = |V(H_2)| = d+1, \ d+1 \leq |V(H_3)| \leq d+2, \ \delta(H) = d \geq 3.$$

Let $H_1, H_2, H_3, H_4$ be disjoint graphs with $x_i, y_i \in V(H_i)$ ($i = 1, 2, 3$) and $V(H_4) = \{z\}$. Form a graph $H_5$ by identifying $x_1, x_2, x_3$ and adding the edges $zy_i$ ($i = 1, 2, 3$) and $y_1 y_2$. We will say that $H \in \Re_2$ if and only if

$$H_5 \subseteq H \subseteq H_5 + \{y_1 y_3, y_2 y_3, z x_1\}, \ |V(H_i)| = 4 \ (i = 1, 2, 3), \ \delta(H) = 3.$$

Let $H_1, H_2, H_3, H_4$ be disjoint graphs with $x_i, y_i \in V(H_i)$ ($i = 1, 2, 3$) and $V(H_4) = \{z\}$. Form a graph $H_5$ by identifying $x_1, x_2, x_3$ and adding the edges $zy_i$ ($i = 1, 2, 3$) and $zx_1, y_1 y_2$. We will say that $H \in \Re_3$ if and only if

$$H_5 \subseteq H \subseteq H_5 + \{y_1 y_3, y_2 y_3\}, \ |V(H_i)| = 5 \ (i = 1, 2, 3), \ \delta(H) = 4.$$

Let $H_1, H_2, H_3, H_4$ be disjoint graphs with $x_i, y_i, z_i \in V(H_i)$ ($i = 1, 2, 3, 4$). Form a graph $H$ by identifying $x_1, x_2, x_3, x_4$ in one vertex and identifying $y_1, y_2, y_3, y_4$ in another vertex, as well as by adding the edges of a triangle between $z_1, z_2, z_3$. We will say that $H \in \Re_4$ if and only if

$$|V(H_i)| = 5 \ (i = 1, 2, 3, 4), \ \delta(H) = 4.$$

The class $\Re$ is the union of the classes of graphs $\Re_1, \Re_2, \Re_3, \Re_4$. It is easy to check that if $G \in \Re$ then $\kappa(G) = 2$, $\tau(G) = 1$ and each longest cycle in $G$ is not a dominating cycle.

**Theorem 1**. Let $G$ be a 1-tough graph. If $\delta \geq \frac{1}{3}(n-2)$ then each longest cycle in $G$ is a dominating cycle if and only if $G \notin \Re$.

The next theorem follows from Theorem 1 immediately, based on the fact that $\kappa(G) = 2$ for each $G \in \Re$.

**Theorem 2**. Let $G$ be a 3-connected 1-tough graph. If $\delta \geq \frac{1}{3}(n-2)$ then each longest cycle in $G$ is a dominating cycle.

To prove Theorems 1 and 2, we need a known lower bound for the circumference, the length of a longest cycle in a graph, concerning the alternative existence of long cycles and Hamilton cycles in 1-tough graphs.

**Theorem D [1].** Every 1-tough graph either has a Hamilton cycle or has a cycle of length at least $2\delta + 2$.

## 2  Notations and preliminaries

The set of vertices of a graph $G$ is denoted by $V(G)$ and the set of edges by $E(G)$. For $S$ a subset of $V(G)$, we denote by $G\backslash S$ the maximum subgraph of $G$ with vertex set $V(G)\backslash S$. We write $G[S]$ for the subgraph of $G$ induced by $S$. For a subgraph $H$ of $G$ we use $G\backslash H$ short for $G\backslash V(H)$. The neighborhood of a vertex $x \in V(G)$ is denoted by $N(x)$. Furthermore, for a subgraph $H$ of $G$ and $x \in V(G)$, we define $N_H(x) = N(x) \cap V(H)$ and $d_H(x) = |N_H(x)|$. If $X \subset V(G)$ then $N(X)$ denotes the set of all vertices of $G\backslash X$ adjacent to vertices in $X$. If $X = \{x_1, x_2, ..., x_r\}$ then $N(X)$ is written as $N(x_1, x_2, ..., x_r)$. Let $s(G)$ denote the number of components of a graph $G$. A graph $G$ is $t$-tough if $|S| \geq ts(G\backslash S)$ for every subset $S$ of the vertex set $V(G)$ with $s(G\backslash S) > 1$. The toughness of $G$, denoted $\tau(G)$, is the maximum value of $t$ for which $G$ is $t$-tough (taking $\tau(K_n) = \infty$ for all $n \geq 1$).

A simple cycle (or just a cycle) $C$ of length $t$ is a sequence $v_1 v_2 ... v_t v_1$ of distinct vertices $v_1, ..., v_t$ with $v_i v_{i+1} \in E(G)$ for each $i \in \{1, ..., t\}$, where $v_{t+1} = v_1$. When $t = 2$, the cycle $C = v_1 v_2 v_1$ on two vertices $v_1, v_2$ coincides with the edge $v_1 v_2$, and when $t = 1$, the cycle $C = v_1$ coincides with the vertex $v_1$. So, all vertices and edges in a graph can be considered as cycles of lengths 1 and 2, respectively. A graph $G$ is hamiltonian if $G$ contains a Hamilton cycle, i.e. a cycle of length $n$. A cycle $C$ in $G$ is dominating if $G\backslash C$ is edgeless.

Paths and cycles in a graph $G$ are considered as subgraphs of $G$. If $Q$ is a path or a cycle, then the length of $Q$, denoted $|Q|$, is $|E(Q)|$. We write $Q$ with a given orientation by $\overrightarrow{Q}$. For $x, y \in V(Q)$, we denote by $x\overrightarrow{Q}y$ the subpath of $Q$ in the chosen direction from $x$ to $y$. For $x \in V(Q)$, we denote the $h$-th successor and the $h$-th predecessor of $x$ on $\overrightarrow{Q}$ by $x^{+h}$ and $x^{-h}$, respectively. We abbreviate $x^{+1}$ and $x^{-1}$ by $x^+$ and $x^-$, respectively.

**Special definitions**. Let $G$ be a graph, $C$ a longest cycle in $G$ and $P = x\overrightarrow{P}y$ a longest path in $G\backslash C$ of length $\overline{p} \geq 0$. Let $\xi_1, \xi_2, ..., \xi_s$ be the elements of $N_C(x) \cup N_C(y)$ occuring on $C$ in a consecutive order. Set

$$I_i = \xi_i \overrightarrow{C} \xi_{i+1}, \ I_i^* = \xi_i^+ \overrightarrow{C} \xi_{i+1}^- \ \ (i = 1, 2, ..., s),$$

where $\xi_{s+1} = \xi_1$.

(1) The segments $I_1, I_2, ..., I_s$ are called elementary segments on $C$ created by $N_C(x) \cup N_C(y)$.

(2) We call a path $L = z\overrightarrow{L}w$ an intermediate path between two distinct elementary segments $I_a$ and $I_b$ if

$$z \in V(I_a^*), \ w \in V(I_b^*), \ V(L) \cap V(C \cup P) = \{z, w\}.$$

(3) Define $\Upsilon(I_{i_1}, I_{i_2}, ..., I_{i_t})$ to be the set of all intermediate paths between elementary segments $I_{i_1}, I_{i_2}, ..., I_{i_t}$.

**Lemma 1.** Let $G$ be a graph, $C$ a longest cycle in $G$ and $P = x\overrightarrow{P}y$ a longest path in $G\backslash C$ of length $\overline{p} \geq 1$. If $|N_C(x)| \geq 2$, $|N_C(y)| \geq 2$ and $N_C(x) \neq N_C(y)$ then

$$|C| \geq \begin{cases} 3\delta + \max\{\sigma_1, \sigma_2\} - 1 \geq 3\delta & \text{if} \quad \overline{p} = 1, \\ \max\{2\overline{p} + 8, 4\delta - 2\overline{p}\} & \text{if} \quad \overline{p} \geq 2, \end{cases}$$

where $\sigma_1 = |N_C(x)\backslash N_C(y)|$ and $\sigma_2 = |N_C(y)\backslash N_C(x)|$.

**Lemma 2.** Let $G$ be a graph, $C$ a longest cycle in $G$ and $P = x\overrightarrow{P}y$ a longest path in $G\backslash C$ of length $\overline{p} \geq 0$. If $N_C(x) = N_C(y)$ and $|N_C(x)| \geq 2$ then for each elementary segments $I_a$ and $I_b$ induced by $N_C(x) \cup N_C(y)$,

(a1) if $L$ is an intermediate path between $I_a$ and $I_b$ then

$$|I_a| + |I_b| \geq 2\overline{p} + 2|L| + 4,$$

(a2) if $\Upsilon(I_a, I_b) \subseteq E(G)$ and $|\Upsilon(I_a, I_b)| = i$ for some $i \in \{1, 2, 3\}$ then

$$|I_a| + |I_b| \geq 2\overline{p} + i + 5,$$

(a3) if $\Upsilon(I_a, I_b) \subseteq E(G)$ and $\Upsilon(I_a, I_b)$ contains two independent intermediate edges then

$$|I_a| + |I_b| \geq 2\overline{p} + 8.$$

## 3 Proofs

**Proof of Lemma 1**. Put

$$A_1 = N_C(x)\backslash N_C(y), \ A_2 = N_C(y)\backslash N_C(x), \ M = N_C(x) \cap N_C(y).$$

By the hypothesis, $N_C(x) \neq N_C(y)$, implying that

$$\max\{|A_1|, |A_2|\} \geq 1.$$

Let $\xi_1, \xi_2, ..., \xi_s$ be the elements of $N_C(x) \cup N_C(y)$ occuring on $C$ in a consecutive order. Put $I_i = \xi_i \overrightarrow{C} \xi_{i+1}$ ($i = 1, 2, ..., s$), where $\xi_{s+1} = \xi_1$. Clearly, $s = |A_1| + |A_2| + |M|$. Since $C$ is extreme, $|I_i| \geq 2$ ($i = 1, 2, ..., s$). Next, if $\{\xi_i, \xi_{i+1}\} \cap M \neq \emptyset$ for some $i \in \{1, 2, ..., s\}$ then $|I_i| \geq \overline{p} + 2$. Further, if either $\xi_i \in A_1$, $\xi_{i+1} \in A_2$ or $\xi_i \in A_2$, $\xi_{i+1} \in A_1$ then again $|I_i| \geq \overline{p} + 2$.

**Case 1**. $\overline{p} = 1$.
**Case 1.1**. $|A_i| \geq 1$ $(i = 1, 2)$.

It follows that among $I_1, I_2, ..., I_s$ there are $|M| + 2$ segments of length at least $\overline{p} + 2$. Observing also that each of the remaining $s - (|M| + 2)$ segments has a length at least 2, we get

$$|C| \geq (\overline{p} + 2)(|M| + 2) + 2(s - |M| - 2)$$

$$= 3(|M| + 2) + 2(|A_1| + |A_2| - 2)$$

$$= 2|A_1| + 2|A_2| + 3|M| + 2.$$

Since $|A_1| = d(x) - |M| - 1$ and $|A_2| = d(y) - |M| - 1$, we have

$$|C| \geq 2d(x) + 2d(y) - |M| - 2 \geq 3\delta + d(x) - |M| - 2.$$

Recalling that $d(x) = |M| + |A_1| + 1$, we obtain

$$|C| \geq 3\delta + |A_1| - 1 = 3\delta + \sigma_1 - 1.$$

Analogously, $|C| \geq 3\delta + \sigma_2 - 1$. So,

$$|C| \geq 3\delta + \max\{\sigma_1, \sigma_2\} - 1 \geq 3\delta.$$

**Case 1.2**. Either $|A_1| \geq 1, |A_2| = 0$ or $|A_1| = 0, |A_2| \geq 1$.

Assume w.l.o.g. that $|A_1| \geq 1$ and $|A_2| = 0$, i.e. $|N_C(y)| = |M| \geq 2$ and $s = |A_1| + |M|$. Hence, among $I_1, I_2, ..., I_s$ there are $|M| + 1$ segments of length at least $\overline{p} + 2 = 3$. Taking into account that each of the remaining $s - (|M| + 1)$ segments has a length at least 2 and $|M| + 1 = d(y)$, we get

$$|C| \geq 3(|M| + 1) + 2(s - |M| - 1) = 3d(y) + 2(|A_1| - 1)$$

$$\geq 3\delta + |A_1| - 1 = 3\delta + \max\{\sigma_1, \sigma_2\} - 1 \geq 3\delta.$$

**Case 2**. $\overline{p} \geq 2$.

We first prove that $|C| \geq 2\overline{p} + 8$. Since $|N_C(x)| \geq 2$ and $|N_C(y)| \geq 2$, there are at least two segments among $I_1, I_2, ..., I_s$ of length at least $\overline{p} + 2$. If $|M| = 0$ then clearly $s \geq 4$ and

$$|C| \geq 2(\overline{p} + 2) + 2(s - 2) \geq 2\overline{p} + 8.$$

Otherwise, since $\max\{|A_1|, |A_2|\} \geq 1$, there are at least three elementary segments of length at least $\overline{p} + 2$, that is

$$|C| \geq 3(\overline{p} + 2) \geq 2\overline{p} + 8.$$

So, in any case, $|C| \geq 2\overline{p} + 8$.

To prove that $|C| \geq 4\delta - 2\overline{p}$, we distinguish two main cases.

**Case 2.1**. $|A_i| \geq 1$ $(i = 1, 2)$.

It follows that among $I_1, I_2, ..., I_s$ there are $|M| + 2$ segments of length at least $\overline{p} + 2$. Further, since each of the remaining $s - (|M| + 2)$ segments has a length at least 2, we get

$$|C| \geq (\overline{p} + 2)(|M| + 2) + 2(s - |M| - 2)$$
$$= (\overline{p} - 2)|M| + (2\overline{p} + 4|M| + 4) + 2(|A_1| + |A_2| - 2)$$
$$\geq 2|A_1| + 2|A_2| + 4|M| + 2\overline{p}.$$

Observing also that

$$|A_1| + |M| + \overline{p} \geq d(x), \quad |A_2| + |M| + \overline{p} \geq d(y),$$

we have

$$2|A_1| + 2|A_2| + 4|M| + 2\overline{p}$$
$$\geq 2d(x) + 2d(y) - 2\overline{p} \geq 4\delta - 2\overline{p},$$

implying that $|C| \geq 4\delta - 2\overline{p}$.

**Case 2.2.** Either $|A_1| \geq 1, |A_2| = 0$ or $|A_1| = 0, |A_2| \geq 1$.

Assume w.l.o.g. that $|A_1| \geq 1$ and $|A_2| = 0$, i.e. $|N_C(y)| = |M| \geq 2$ and $s = |A_1| + |M|$. It follows that among $I_1, I_2, ..., I_s$ there are $|M| + 1$ segments of length at least $\overline{p} + 2$. Observing also that $|M| + \overline{p} \geq d(y) \geq \delta$, i.e. $2\overline{p} + 4|M| \geq 4\delta - 2\overline{p}$, we get

$$|C| \geq (\overline{p} + 2)(|M| + 1) \geq (\overline{p} - 2)(|M| - 1) + 2\overline{p} + 4|M|$$
$$\geq 2\overline{p} + 4|M| \geq 4\delta - 2\overline{p}. \quad \blacksquare$$

**Proof of Lemma 2.** Let $\xi_1, \xi_2, ..., \xi_s$ be the elements of $N_C(x)$ occuring on $C$ in a consecutive order. Put $I_i = \xi_i \overrightarrow{C} \xi_{i+1}$ $(i = 1, 2, ..., s)$, where $\xi_{s+1} = \xi_1$. To prove (a1), let $L = z\overrightarrow{L}w$ be an intermediate path between elementary segments $I_a$ and $I_b$ with $z \in V(I_a^*)$ and $w \in V(I_b^*)$. Put

$$|\xi_a \overrightarrow{C} z| = d_1, \ |z\overrightarrow{C}\xi_{a+1}| = d_2, \ |\xi_b \overrightarrow{C} w| = d_3, \ |w\overrightarrow{C}\xi_{b+1}| = d_4,$$
$$C' = \xi_a x \overrightarrow{P} y \xi_b \overleftarrow{C} z \overrightarrow{L} w \overrightarrow{C} \xi_a.$$

Clearly,
$$|C'| = |C| - d_1 - d_3 + |L| + |P| + 2.$$

Since $C$ is extreme, we have $|C| \geq |C'|$, implying that $d_1 + d_3 \geq \overline{p} + |L| + 2$. By a symmetric argument, $d_2 + d_4 \geq \overline{p} + |L| + 2$. Hence

$$|I_a| + |I_b| = \sum_{i=1}^{4} d_i \geq 2\overline{p} + 2|L| + 4.$$

The proof of $(a1)$ is complete. To proof $(a2)$ and $(a3)$, let $\Upsilon(I_a, I_b) \subseteq E(G)$ and $|\Upsilon(I_a, I_b)| = i$ for some $i \in \{1, 2, 3\}$.

**Case 1.** $i = 1$.

It follows that $\Upsilon(I_a, I_b)$ consists of a unique intermediate edge $L = zw$. By $(a1)$,
$$|I_a| + |I_b| \geq 2\overline{p} + 2|L| + 4 = 2\overline{p} + 6.$$

**Case 2.** $i = 2$.

It follows that $\Upsilon(I_a, I_b)$ consists of two edges $e_1, e_2$. Put $e_1 = z_1 w_1$ and $e_2 = z_2 w_2$, where $\{z_1, z_2\} \subseteq V(I_a^*)$ and $\{w_1, w_2\} \subseteq V(I_b^*)$.

**Case 2.1.** $z_1 \neq z_2$ and $w_1 \neq w_2$.

Assume w.l.o.g. that $z_1$ and $z_2$ occur in this order on $I_a$.

**Case 2.1.1.** $w_2$ and $w_1$ occur in this order on $I_b$.

Put
$$|\xi_a \overrightarrow{C} z_1| = d_1,\ |z_1 \overrightarrow{C} z_2| = d_2,\ |z_2 \overrightarrow{C} \xi_{a+1}| = d_3,$$
$$|\xi_b \overrightarrow{C} w_2| = d_4,\ |w_2 \overrightarrow{C} w_1| = d_5,\ |w_1 \overrightarrow{C} \xi_{b+1}| = d_6,$$
$$C' = \xi_a \overrightarrow{C} z_1 w_1 \overleftarrow{C} w_2 z_2 \overrightarrow{C} \xi_b x \overrightarrow{P} y \xi_{b+1} \overrightarrow{C} \xi_a.$$

Clearly,
$$|C'| = |C| - d_2 - d_4 - d_6 + |\{e_1\}| + |\{e_2\}| + |P| + 2$$
$$= |C| - d_2 - d_4 - d_6 + \overline{p} + 4.$$

Since $C$ is extreme, $|C| \geq |C'|$, implying that $d_2 + d_4 + d_6 \geq \overline{p} + 4$. By a symmetric argument, $d_1 + d_3 + d_5 \geq \overline{p} + 4$. Hence
$$|I_a| + |I_b| = \sum_{i=1}^{6} d_i \geq 2\overline{p} + 8.$$

**Case 2.1.2.** $w_1$ and $w_2$ occur in this order on $I_b$.

Putting
$$C' = \xi_a \overrightarrow{C} z_1 w_1 \overrightarrow{C} w_2 z_2 \overrightarrow{C} \xi_b x \overrightarrow{P} y \xi_{b+1} \overrightarrow{C} \xi_a,$$
we can argue as in Case 2.1.1.

**Case 2.2.** Either $z_1 = z_2$, $w_1 \neq w_2$ or $z_1 \neq z_2$, $w_1 = w_2$.

Assume w.l.o.g. that $z_1 \neq z_2$, $w_1 = w_2$ and $z_1, z_2$ occur in this order on $I_a$. Put
$$|\xi_a \overrightarrow{C} z_1| = d_1,\ |z_1 \overrightarrow{C} z_2| = d_2,\ |z_2 \overrightarrow{C} \xi_{a+1}| = d_3,$$
$$|\xi_b \overrightarrow{C} w_1| = d_4,\ |w_1 \overrightarrow{C} \xi_{b+1}| = d_5,$$
$$C' = \xi_a x \overrightarrow{P} y \xi_b \overleftarrow{C} z_1 w_1 \overrightarrow{C} \xi_a,$$
$$C'' = \xi_a \overrightarrow{C} z_2 w_1 \overleftarrow{C} \xi_{a+1} x \overrightarrow{P} y \xi_{b+1} \overrightarrow{C} \xi_a.$$

Clearly,
$$|C'| = |C| - d_1 - d_4 + |\{e_1\}| + |P| + 2 = |C| - d_1 - d_4 + \overline{p} + 3,$$
$$|C''| = |C| - d_3 - d_5 + |\{e_2\}| + |P| + 2 = |C| - d_3 - d_5 + \overline{p} + 3.$$
Since $C$ is extreme, $|C| \geq |C'|$ and $|C| \geq |C''|$, implying that
$$d_1 + d_4 \geq \overline{p} + 3, \ d_3 + d_5 \geq \overline{p} + 3.$$
Hence,
$$|I_a| + |I_b| = \sum_{i=1}^{5} d_i \geq d_1 + d_3 + d_4 + d_5 + 1 \geq 2\overline{p} + 7.$$

**Case 3.** $i = 3$.

It follows that $\Upsilon(I_a, I_b)$ consists of three edges $e_1, e_2, e_3$. Let $e_i = z_i w_i$ ($i = 1, 2, 3$), where $\{z_1, z_2, z_3\} \subseteq V(I_a^*)$ and $\{w_1, w_2, w_3\} \subseteq V(I_b^*)$. If there are two independent edges among $e_1, e_2, e_3$ then we can argue as in Case 2.1. Otherwise, we can assume w.l.o.g. that $w_1 = w_2 = w_3$ and $z_1, z_2, z_3$ occur in this order on $I_a$. Put
$$|\xi_a \overrightarrow{C} z_1| = d_1, \ |z_1 \overrightarrow{C} z_2| = d_2, \ |z_2 \overrightarrow{C} z_3| = d_3,$$
$$|z_3 \overrightarrow{C} \xi_{a+1}| = d_4, \ |\xi_b \overrightarrow{C} w_1| = d_5, \ |w_1 \overrightarrow{C} \xi_{b+1}| = d_6,$$
$$C' = \xi_a x \overrightarrow{P} y \xi_b \overleftarrow{C} z_1 w_1 \overrightarrow{C} \xi_a,$$
$$C'' = \xi_a \overrightarrow{C} z_3 w_1 \overleftarrow{C} \xi_{a+1} x \overrightarrow{P} y \xi_{b+1} \overrightarrow{C} \xi_a.$$
Clearly,
$$|C'| = |C| - d_1 - d_5 + |\{e_1\}| + \overline{p} + 2,$$
$$|C''| = |C| - d_4 - d_6 + |\{e_3\}| + \overline{p} + 2.$$
Since $C$ is extreme, we have $|C| \geq |C'|$ and $|C| \geq |C''|$, implying that
$$d_1 + d_5 \geq \overline{p} + 3, \ d_4 + d_6 \geq \overline{p} + 3.$$
Hence,
$$|I_a| + |I_b| = \sum_{i=1}^{6} d_i \geq d_1 + d_4 + d_5 + d_6 + 2 \geq 2\overline{p} + 8. \ \blacksquare$$

**Proof of Theorem 1.** Let $C$ be a longest cycle in $G$ and $P = x_1 \overrightarrow{P} x_2$ a longest path in $G \backslash C$ of length $\overline{p}$. If $|V(P)| \leq 1$ then $C$ is a dominating cycle and we are done. Let $|V(P)| \geq 2$, that is $\overline{p} \geq 1$. Since $\tau \geq 1$, $G$ is 2-connected. By Theorem D, $|C| \geq 2\delta + 2$. On the other hand, by the hypothesis, $|C| + \overline{p} + 1 \leq n \leq 3\delta + 2$, implying that
$$2\delta + 2 \leq |C| \leq 3\delta - \overline{p} + 1. \tag{1}$$

Let $\xi_1, \xi_2, ..., \xi_s$ be the elements of $N_C(x_1) \cup N_C(x_2)$ occuring on $C$ in a consecutive order. Put

$$I_i = \xi_i \overrightarrow{C} \xi_{i+1}, \ I_i^* = \xi_i^+ \overrightarrow{C} \xi_{i+1}^- \ (i = 1, 2, ..., s),$$

where $\xi_{s+1} = \xi_1$. Assume first that $\delta = 2$. By (1), $\overline{p} = 1$, $|C| = 6$ and $n = 8$. Then it is easy to see that $s = 2$, $|I_1| = |I_2| = 3$ and $\Upsilon(I_1, I_2) = \emptyset$. This means that $G \backslash \{\xi_1, \xi_2\}$ has three connected components, that is $\tau < 1$, contradicting the hypothesis. Now let $\delta \geq 3$. By (1), $\overline{p} \leq \delta - 1$.

**Case 1.** $\overline{p} = 1$.
By (1),
$$2\delta + 2 \leq |C| \leq 3\delta. \tag{2}$$
Since $\delta \geq 3$, we have $|N_C(x_i)| \geq \delta - \overline{p} = \delta - 1 \geq 2 \ (i = 1, 2)$.

**Case 1.1.** $N_C(x_1) \neq N_C(x_2)$.
It follows that $\max\{\sigma_1, \sigma_2\} \geq 1$, where

$$\sigma_1 = |N_C(x_1) \backslash N_C(x_2)|, \quad \sigma_2 = |N_C(x_2) \backslash N_C(x_1)|.$$

By Lemma 1, $|C| \geq 3\delta$. Recalling (2), we get $|C| = 3\delta$. If $\max\{\sigma_1, \sigma_2\} \geq 2$ then by Lemma 1, $|C| \geq 3\delta + 1$, contradicting (2). Let $\max\{\sigma_1, \sigma_2\} = 1$. This implies $s \geq \delta$ and $|I_i| \geq 3 \ (i = 1, 2, ..., s)$. If $s \geq \delta + 1$ then $|C| \geq 3s \geq 3\delta + 3$, again contradicting (2). Hence $s = \delta$, which yields $|I_i| = 3 \ (i = 1, 2, ..., s)$. Assume that $\Upsilon(I_1, I_2, ..., I_s) \neq \emptyset$, that is $\Upsilon(I_a, I_b) \neq \emptyset$ for some distinct $a, b \in \{1, 2, ..., s\}$. By the definition, there is an intermediate path $L = y\overrightarrow{L}z$ between $I_a$ and $I_b$. Put

$$|\xi_a \overrightarrow{C} y| = d_1, \ |y \overrightarrow{C} \xi_{a+1}| = d_2, \ |\xi_b \overrightarrow{C} z| = d_3, \ |z \overrightarrow{C} \xi_{b+1}| = d_4.$$

Let $Q_1 = \xi_a \overrightarrow{Q}_1 \xi_b$ and $Q_2 = \xi_{a+1} \overrightarrow{Q}_2 \xi_{b+1}$ be two longest paths in $G$ with

$$V(Q_1) \cap V(C) = \{\xi_a, \xi_b\}, \quad V(Q_2) \cap V(C) = \{\xi_{a+1}, \xi_{b+1}\}.$$

Since $\max\{\sigma_1, \sigma_2\} = 1$, it is easy to see that $|Q_i| \geq 3 \ (i = 1, 2)$. Consequently,

$$|C| \geq |\xi_a Q_1 \xi_b \overleftarrow{C} y \overrightarrow{L} z \overrightarrow{C} \xi_a| = |C| - d_1 - d_3 + |Q_1| + |L|,$$

implying that $d_1 + d_3 \geq |Q_1| + |L| \geq 4$. By a symmetric argument, $d_2 + d_4 \geq 4$. By summing, we get $|I_a| + |I_b| = \sum_{i=1}^{4} d_i \geq 8$, contradicting the fact that $|I_a| = |I_b| = 3$. So, $\Upsilon(I_1, I_2, ..., I_s) = \emptyset$. This means that $G \backslash \{\xi_1, \xi_2, ..., \xi_s\}$ has at least $s + 1$ connected components, contradicting the fact that $\tau \geq 1$.

**Case 1.2.** $N_C(x_1) = N_C(x_2)$.
Clearly, $s = |N_C(x_1)| \geq \delta - \overline{p} = \delta - 1 \geq 2$ and $|I_i| \geq \overline{p} + 2 = 3 \ (i = 1, 2, ..., s)$. If $s \geq \delta$ then $|C| \geq 3s \geq 3\delta$ and we can argue as in Case 1.1. Let $s = \delta - 1$.
The following claim can be derived from (2) easily.

**Claim 1**. (1) $|I_i| + |I_j| \leq 9$ for each distinct $i, j \in \{1, 2, ..., s\}$.

(2) If $|I_a| + |I_b| = 9$ for some distinct $a, b \in \{1, 2, ..., s\}$ then $|I_i| = 3$ for each $i \in \{1, 2, ..., s\}\setminus\{a, b\}$.

(3) If $|I_a| = 6$ for some $a \in \{1, 2, ..., s\}$ then $|I_i| = 3$ for each $i \in \{1, 2, ..., s\}\setminus\{a\}$.

(4) There are at most three elementary segments of length at least 4.

(5) If $|I_a| \geq 4$, $|I_b| \geq 4$, $|I_c| \geq 4$ for some distinct $a, b, c \in \{1, 2, ..., s\}$ then $|I_a| = |I_b| = |I_c| = 4$ and $|I_i| = 3$ for each $i \in \{1, 2, ..., s\}\setminus\{a, b, c\}$.

(6) If $|I_a| \geq 4$ and $|I_b| \geq 4$ for some distinct $a, b \in \{1, 2, ..., s\}$ then $|I_i| \leq 4$ for each $i \in \{1, 2, ..., s\}$.

If $\Upsilon(I_1, I_2, ..., I_s) = \emptyset$ then $G\setminus\{\xi_1, \xi_2, ..., \xi_s\}$ has at least $s + 1$ connected component, contradicting the fact that $\tau \geq 1$. Otherwise, $\Upsilon(I_a, I_b) \neq \emptyset$ for some distinct $a, b \in \{1, 2, ..., s\}$. By definition, there is an intermediate path $L$ between $I_a$ and $I_b$. If $|L| \geq 2$ then by Lemma 2(a1),

$$|I_a| + |I_b| \geq 2\overline{p} + 2|L| + 4 \geq 10,$$

contradicting Claim 1(1). Otherwise, $|L| = 1$ and therefore,

$$\Upsilon(I_1, I_2, ..., I_s) \subseteq E(G).$$

By Lemma 2(a1), $|I_a| + |I_b| \geq 2\overline{p} + 6 = 8$. Combining this with Claim 1(1), we have

$$8 \leq |I_a| + |I_b| \leq 9.$$

Let $L = yz$, where $y \in V(I_a^*)$ and $z \in V(I_b^*)$. Put

$$C_1 = \xi_a x_1 x_2 \xi_b \overleftarrow{C} yz \overrightarrow{C} \xi_a,$$

$$C_2 = \xi_a \overrightarrow{C} yz \overleftarrow{C} \xi_{a+1} x_1 x_2 \xi_{b+1} \overrightarrow{C} \xi_a.$$

**Case 1.2.1.** $|I_a| + |I_b| = 8$.

If $|\Upsilon(I_a, I_b)| \geq 2$ then by Lemma 2(a2), $|I_a| + |I_b| \geq 9$, contradicting the hypothesis. Hence,

$$\Upsilon(I_a, I_b) = \{yz\}. \qquad (3)$$

Since $|I_i| \geq 3$ ($i = 1, 2, ..., s$), we can assume w.l.o.g. that either $|I_a| = 3, |I_b| = 5$ or $|I_a| = |I_b| = 4$.

**Case 1.2.1.1.** $|I_a| = 3$ and $|I_b| = 5$.

Put $I_a = \xi_a w_1 w_2 \xi_{a+1}$ and $I_b = \xi_b w_3 w_4 w_5 w_6 \xi_{b+1}$. Assume w.l.o.g. that $y = w_2$. If $z = w_3$ then $|C_1| > |C|$, a contradiction. Further, if $z \in \{w_5, w_6\}$ then $|C_2| > |C|$, again a contradiction. Hence, $z = w_4$. Further, if $\Upsilon(I_a, I_c) \neq \emptyset$ for some $c \in \{1, 2, ..., s\}\setminus\{a, b\}$ then by Lemma 2(a1), $|I_a| + |I_c| \geq 2\overline{p} + 6 = 8$, implying that $|I_c| \geq 5$. But then $|I_b| + |I_c| \geq 10$, contradicting Claim 1(1). Hence $\Upsilon(I_a, I_i) = \emptyset$ for each $i \in \{1, 2, ..., s\}\setminus\{a, b\}$. Combining this with (3), we get $N(w_1) \subseteq \{\xi_1, \xi_2, ..., \xi_s, w_2\}$. Furthermore, if $w_1 \xi_{a+1} \in E(G)$ then

$$\xi_a x_1 x_2 \xi_b \overleftarrow{C} \xi_{a+1} w_1 w_2 w_4 \overrightarrow{C} \xi_a$$

is longer than $C$, a contradiction. Thus $|N(w_1)| \leq s = \delta - 1$, contradicting the fact that $|N(w_1)| = d(w_1) \geq \delta$.

**Case 1.2.1.2.** $|I_a| = |I_b| = 4$.
Put $I_a = \xi_a w_1 w_2 w_3 \xi_{a+1}$ and $I_b = \xi_b w_4 w_5 w_6 \xi_{b+1}$.

**Case 1.2.1.2.1.** $y \in \{w_1, w_3\}$.
Assume w.l.o.g. that $y = w_3$. If $z \in \{w_5, w_6\}$ than $|C_2| > |C|$, a contradiction. Hence $z = w_4$. If $u\xi_{a+1} \in E(G)$ for some $u \in \{w_1, w_2\}$, then

$$|\xi_a x_1 x_2 \xi_b \overleftarrow{C} \xi_{a+1} u \overrightarrow{C} w_3 w_4 \overrightarrow{C} \xi_a| \geq |C| + 1,$$

a contradiction. Next, if $u\xi_b \in E(G)$ then

$$|\xi_a x_1 x_2 \xi_{a+1} \overrightarrow{C} \xi_b u \overrightarrow{C} w_3 w_4 \overrightarrow{C} \xi_a| \geq |C| + 1,$$

a contradiction. Thus, if $u \in \{w_1, w_2\}$ then $N(u) \cap \{\xi_{a+1}, \xi_b\} = \emptyset$. By a symmetric argument, if $v \in \{w_5, w_6\}$ then $N(v) \cap \{\xi_{a+1}, \xi_b\} = \emptyset$. Combining these relations with (3), we get

$$N(u) \cap (\{\xi_{a+1}, \xi_b\} \cup V(I_b^*)) = \emptyset \quad \text{for each} \quad u \in \{w_1, w_2\}, \tag{4}$$

$$N(u) \cap (\{\xi_{a+1}, \xi_b\} \cup V(I_a^*)) = \emptyset \quad \text{for each} \quad u \in \{w_5, w_6\}. \tag{5}$$

**Case 1.2.1.2.1.1.** $\xi_{a+1} \neq \xi_b$.
If $N(u) \cap V(I_i^*) = \emptyset$ for some $u \in \{w_1, w_2\}$ and $i \in \{1, 2, ..., s\}\backslash\{a\}$ then by (4),

$$N(u) \subseteq \{\xi_1, \xi_2, ..., \xi_s, w_1, w_2, w_3\}\backslash\{\xi_{a+1}, \xi_b, u\},$$

contradicting the fact that $|N(u)| \geq \delta = s+1$. Otherwise, by (3), $w_1 v \in E(G)$, where $v \in V(I_c^*)$ for some $c \in \{1, 2, ..., s\}\backslash\{a, b\}$, and $w_2 u \in E(G)$, where $u \in V(I_d^*)$ for some $d \in \{1, 2, ..., s\}\backslash\{a, b\}$. By Lemma 2(a1), $|I_a| + |I_c| \geq 2\overline{p} + 6 = 8$, implying that $|I_c| \geq 4$. By Claim 1(5), $|I_c| = 4$. By the same reason, $|I_d| = 4$. If $d = c$ then $|\Upsilon(I_a, I_c)| = 2$ and by Lemma 2(a2), $|I_a| + |I_c| \geq 2\overline{p} + 7 = 9$, contradicting the fact that $|I_a| = |I_c| = 4$. Otherwise, there are at least four elementary segments of length at least 4, contradicting Claim 1(4).

**Case 1.2.1.2.1.2.** $\xi_{a+1} = \xi_b$.
Assume w.l.o.g. that $a = 1$ and $b = 2$.

**Case 1.2.1.2.1.2.1.** $s = 2$.
It follows that $\delta = s + 1 = 3$, $|C| = |I_1| + |I_2| = 8$ and $10 \leq n \leq 11$. By (4) and (5), $N(\xi_2) \cap \{w_1, w_2, w_5, w_6\} = \emptyset$. If $n = 10$ then $V(G) = V(C \cup P)$ and therefore,

$$N(x_1, x_2) = \{\xi_1, \xi_2\}, \ N(w_1, w_2) = \{\xi_1, w_3\}, \ N(w_5, w_6) = \{\xi_1, w_4\},$$

implying that $G \in \Re_1$. Now let $n = 11$, that is there is a vertex $x_3 \in V(G)\backslash V(C \cup P)$. Clearly, $N(x_3) \subseteq V(C)$. Assume first that either $N(x_3) \subseteq$ $V(I_1)$ or $N(x_3) \subseteq V(I_2)$, say $N(x_3) \subseteq V(I_1)$. Since $C$ is extreme, we have $N(x_3) = \{\xi_1, \xi_2, w_2\}$. Then

$$\xi_1 x_1 x_2 \xi_2 x_3 w_2 w_3 w_4 \overrightarrow{C} \xi_1$$

is longer than $C$, a contradiction. Now let $N(x_3) \not\subseteq V(I_i)$ ($i = 1, 2$), that is there is an intermediate path of length 2 between $I_1$ and $I_2$. By Lemma 2(a1), $|I_1| + |I_2| \geq 2\overline{p} + 2|L| + 4 = 10$, a contradiction.

**Case 1.2.1.2.1.2.2.** $s = 3$.

It follows that $\delta = s + 1 = 4$. Assume first that $|I_3| = 3$. Put $I_3 = \xi_3 w_7 w_8 \xi_1$. If $w_8 \xi_2 \in E(G)$ then

$$\xi_1 \overrightarrow{C} w_3 w_4 \overrightarrow{C} w_8 \xi_2 x_2 x_1 \xi_1$$

is longer than $C$, a contradiction. Let $w_8 \xi_2 \notin E(G)$. Observing that by Lemma 2(a1), $\Upsilon(I_1, I_3) = \Upsilon(I_2, I_3) = \emptyset$, we conclude that $N(w_8) \subseteq \{\xi_1, \xi_3, w_7\}$, which contradicts the fact that $|N(w_8)| \geq \delta = 4$. Now assume that $|I_3| \geq 4$. By Claim 2(5), $|I_3| = 4$, implying that $|C| = 12$, $n = 14$ and $V(G) = V(C \cup P)$. Put $I_3 = \xi_3 w_7 w_8 w_9 \xi_1$. If $\xi_2 v \in E(G)$ for some $v \in \{w_8, w_9\}$ then

$$\xi_1 \overrightarrow{C} w_3 w_4 \overrightarrow{C} v \xi_2 x_2 x_1 \xi_1$$

is longer than $C$, a contradiction. Hence $w_8 \xi_2, w_9 \xi_2 \notin E(G)$. By a symmetric argument, $w_7 \xi_2 \notin E(G)$. Recalling also (4), we conclude that $N(\xi_2) \subseteq \{\xi_1, \xi_3, w_3, w_4, x_1, x_2\}$. If $\Upsilon(I_1, I_2, I_3) = \{w_3 w_4\}$ then

$$N(x_1, x_2) = \{\xi_1, \xi_2, \xi_3\}, \; N(w_1, w_2) = \{\xi_1, w_3\},$$

$$N(w_5, w_6) = \{\xi_3, w_4\}, \; N(w_7, w_8, w_9) = \{\xi_1, \xi_3\},$$

implying that $G \in \Re_4$. Otherwise, by (3), there is an intermediate edge $uv$ such that $u \in V(I_1^*) \cup V(I_2^*)$ and $v \in V(I_3^*)$. Assume w.l.o.g. that $u \in V(I_1^*)$. If $u = w_3$ then as in Case 1.2.1.2.1, $v = w_4$, contradicting the fact that $v \in V(I_3^*)$. If $u = w_1$ then, again as in Case 1.2.1.2.1, $v = w_9$ and therefore,

$$|\xi_1 w_1 w_9 \overleftarrow{C} w_4 w_3 \xi_2 x_2 x_1 \xi_1| \geq |C| + 1,$$

a contradiction. Finally, let $u = w_2$. If $v \in \{w_8, w_9\}$ then

$$|\xi_1 w_1 w_2 v \overleftarrow{C} w_4 w_3 \xi_2 x_2 x_1 \xi_1| \geq |C| + 1,$$

a contradiction. If $v = w_7$ then

$$|\xi_1 x_1 x_2 \xi_3 \overleftarrow{C} w_2 w_7 \overrightarrow{C} \xi_1| = |C| + 1,$$

again a contradiction.

**Case 1.2.1.2.1.2.3.** $s \geq 4$.

By Claim 1(4), $|I_c| = 3$ for some $c \in \{3, 4, ..., s\}$. Put $I_c = \xi_c w_7 w_8 \xi_{c+1}$. By Claim 1(6), $|I_i| \leq 4$ ($i = 1, 2, ..., s$). If $w_8 \xi_2 \in E(G)$ then

$$\xi_1 \overrightarrow{C} w_3 w_4 \overrightarrow{C} w_8 \xi_2 x_2 x_1 \xi_{c+1} \overrightarrow{C} \xi_1$$

is longer than $C$, a contradiction. Let $w_8 \xi_2 \notin E(G)$. Observing also that by Lemma 2(a1), $\Upsilon(I_c, I_i) = \emptyset$ for each $i \in \{1, 2, ..., s\} \backslash \{c\}$, we conclude that $N(w_8) \subseteq \{\xi_1, \xi_2, ..., \xi_s, w_7\} \backslash \{\xi_2\}$, which contradicts the fact that $|N(w_8)| \geq \delta = s + 1$.

**Case 1.2.1.2.2.** $y = w_2$.

If $z = w_4$ then $|C_1| > |C|$, a contradiction. Next, if $z = w_6$ then $|C_2| > |C|$, a contradiction. Hence $z = w_5$. We have $w_1 w_3 \notin E(G)$ since otherwise

$$\xi_a x_1 x_2 \xi_b \overleftarrow{C} w_3 w_1 w_2 w_5 \overrightarrow{C} \xi_a$$

is longer than $C$, a contradiction. Further, if $\Upsilon(I_1, I_2, ..., I_s) = \{w_2 w_5\}$ then $G \backslash \{\xi_1, \xi_2, ..., \xi_s, w_2\}$ has at least $s + 2$ connected components, contradicting the fact that $\tau \geq 1$. Now let $\Upsilon(I_1, I_2, ..., I_s) \neq \{w_2 w_5\}$. By (3), $\Upsilon(I_c, I_d) \neq \emptyset$, where $\{c, d\} \neq \{a, b\}$. If $\{c, d\} \cap \{a, b\} = \emptyset$ then by Lemma 2(a1), $|I_c| + |I_d| \geq 8$. Assume w.l.o.g. that $|I_c| \geq 4$. Then by Claim 1(5), $|I_c| = 4$, that is $|I_d| \geq 4$, contradicting Claim 1(4). Now let $\{c, d\} \cap \{a, b\} \neq \emptyset$ and assume w.l.o.g. that $\Upsilon(I_a, I_c) \neq \emptyset$ for some $c \in \{1, 2, ..., s\} \backslash \{a, b\}$. By Lemma 2(a1) and Claim 1(5), $|I_c| = 4$. Put $I_c = \xi_c w_7 w_8 w_9 \xi_{c+1}$. By the definition, there is an intermediate edge $uv$ such that $u \in V(I_a^*)$ and $v \in V(I_c^*)$. If $u \in \{w_1, w_3\}$ then we can argue as in Case 1.2.1.2.1. Let $u = w_2$, implying that $v = w_8$. Recalling that $w_1 w_3 \notin E(G)$ and using symmetric arguments, we can state that $w_9 w_7 \notin E(G)$. Since $|I_a| = |I_b| = |I_c| = 4$, by Claim 1(5) and Lemma 2(a1), $\Upsilon(I_1, I_2, ..., I_s) \subseteq \Upsilon(I_a, I_b, I_c)$. By the same reason, if $\Upsilon(I_a, I_b, I_c) \neq \{w_2 w_5, w_2 w_8\}$, then $\Upsilon(I_a, I_b, I_c) = \{w_2 w_5, w_2 w_8, w_5 w_8\}$. Actually, if $\Upsilon(I_a, I_b, I_c) \neq \{w_2 w_5, w_2 w_8\}$ then $G \backslash \{\xi_1, \xi_2, ..., \xi_s, w_2\}$ has at least $s + 2$ connected components, contradicting the fact that $\tau \geq 1$. Finally, if $\Upsilon(I_a, I_b, I_c) = \{w_2 w_5, w_2 w_8, w_5 w_8\}$ then $G \backslash \{\xi_1, \xi_2, ..., \xi_s, w_2, w_8\}$ has at least $s + 3$ connected components, again contradicting the fact that $\tau \geq 1$.

**Case 1.2.2.** $|I_a| + |I_b| = 9$.

Since $|I_i| \geq 3$ ($i = 1, 2, ..., s$), we can assume w.l.o.g. that either $|I_a| = 3$, $|I_b| = 6$ or $|I_a| = 4$, $|I_b| = 5$.

**Case 1.2.2.1.** $|I_a| = 3$ and $|I_b| = 6$.

By Claim 1(3), $|I_i| = 3$ for each $i \in \{1, 2, ..., s\} \backslash \{b\}$. Put

$$I_a = \xi_a w_1 w_2 \xi_{a+1}, \quad I_b = \xi_b w_3 w_4 w_5 w_6 w_7 \xi_{b+1}.$$

Since $|I_a| = 3$, we can assume w.l.o.g. that $y = w_2$. If $z = w_3$ then $|C_1| > |C|$, a contradiction. If $z \in \{w_6, w_7\}$ then $|C_2| > |C|$, a contradiction. So, $z \in \{w_4, w_5\}$.

**Case 1.2.2.1.1**. $z = w_4$.
If $w_1w_4 \in E(G)$ then

$$|\xi_a x_1 x_2 \xi_b \overleftarrow{C} w_1 w_4 \overrightarrow{C} \xi_a| \geq |C| + 1,$$

a contradiction. Let $w_1w_4 \notin E(G)$. Next, if $N(w_1) \cap \{w_3, w_5, w_6, w_7\} \neq \emptyset$ then there are two independent intermediate edges between $I_a$ and $I_b$ and by Lemma 2(a3), $|I_a| + |I_b| \geq 2\overline{p} + 8 = 10$, contradicting the hypothesis. Let $N(w_1) \cap \{w_3, w_5, w_6, w_7\} = \emptyset$. Further, if $w_1\xi_{a+1} \in E(G)$ then

$$\xi_a x_1 x_2 \xi_b \overleftarrow{C} \xi_{a+1} w_1 w_2 w_4 \overrightarrow{C} \xi_a$$

is longer than $C$, a contradiction. Finally, by Lemma 2(a1) and Claim 1(3), $N(w_1) \cap V(I_i^*) = \emptyset$ for each $i \in \{1, 2, ..., s\} \backslash \{a, b\}$. So,

$$N(w_1) \subseteq \{\xi_1, \xi_2, ..., \xi_s, w_2\} \backslash \{\xi_{a+1}\},$$

contradicting the fact that $|N(w_1)| \geq \delta = s + 1$.

**Case 1.2.2.1.2**. $z = w_5$.
If $w_2w_4 \in E(G)$ then we can argue as in Case 1.2.2.1.1. Let $w_2w_4 \notin E(G)$. Hence, $N(w_2) \cap V(I_b^*) = \{w_5\}$. By a symmetric argument, $N(w_1) \cap V(I_b^*) = \{w_5\}$. Since $|I_i| = |I_a| = 3$ for each $i \in \{1, 2, ..., s\} \backslash \{b\}$, we conclude that $w_5$ is a common vertex for all edges in $\Upsilon(I_1, I_2, ..., I_s)$. Now we consider all possible intermediate connections between $\{w_3, w_4\}$ and $\{w_6, w_7\}$. If $w_3w_6 \in E(G)$ then

$$\xi_a x_1 x_2 \xi_b \overleftarrow{C} w_2 w_5 w_4 w_3 w_6 \overrightarrow{C} \xi_a$$

is longer than $C$, a contradiction. Next, if $w_3w_7 \in E(G)$ then

$$\xi_a w_1 w_2 w_5 w_6 w_7 w_3 \overleftarrow{C} \xi_{a+1} x_1 x_2 \xi_{b+1} \overrightarrow{C} \xi_a$$

is longer than $C$, a contradiction. Finally, if $w_4v \in E(G)$ for some $v \in \{w_6, w_7\}$ then

$$\xi_a w_1 w_2 w_5 \overrightarrow{C} v w_4 \overleftarrow{C} \xi_{a+1} x_1 x_2 \xi_{b+1} \overrightarrow{C} \xi_a$$

is longer than $C$, a contradiction. So, there are no edges between $\{w_3, w_4\}$ and $\{w_6, w_7\}$. This means that $G \backslash \{\xi_1, \xi_2, ..., \xi_s, w_5\}$ has at least $s + 2$ connected components, contradicting the fact that $\tau \geq 1$.

**Case 1.2.2.2**. $|I_a| = 4$ and $|I_b| = 5$.
By Claim 1(2) and Lemma 2(a1), $|I_i| = 3$ and $\Upsilon(I_a, I_i) = \emptyset$ for each $i \in \{1, 2, ..., s\} \backslash \{a, b\}$. If $|\Upsilon(I_a, I_b)| \geq 3$ then by Lemma 2(a2), $|I_a| + |I_b| \geq 2\overline{p} + 8 = 10$, contradicting the hypothesis. Let $1 \leq |\Upsilon(I_a, I_b)| \leq 2$.

**Case 1.2.2.2.1**. $y \in \{w_1, w_3\}$.
Assume w.l.o.g. that $y = w_3$. If $z \in \{w_6, w_7\}$ then $|C_2| > |C|$, a contradiction. Let $z \in \{w_4, w_5\}$.

**Claim 2.** $N(w_1) \subseteq \{\xi_1, \xi_2, ..., \xi_s, w_2, w_3\} \backslash \{\xi_{a+1}, \xi_b\}$.

**Proof.** If $w_1 \xi_{a+1} \in E(G)$ then

$$\xi_a x_1 x_2 \xi_b \overleftarrow{C} \xi_{a+1} w_1 w_2 w_3 z \overrightarrow{C} \xi_a$$

is longer than $C$, a contradiction. Let $w_1 \xi_{a+1} \notin E(G)$. Analogously, $w_1 \xi_b \notin E(G)$. Next, if $w_1 v \in E(G)$ for some $v \in \{w_4, w_5\}$ then

$$\xi_a x_1 x_2 \xi_b \overleftarrow{C} w_1 v \overrightarrow{C} \xi_a$$

is longer than $C$, a contradiction. Hence, $w_1 v \notin E(G)$ for each $v \in \{w_4, w_5\}$. Further, if $w_1 v \in E(G)$ for some $v \in \{w_6, w_7\}$ then $yz$ and $w_1 v$ are independent intermediate edges between $I_a$ and $I_b$. By Lemma 2(a3), $|I_a| + |I_b| \geq 2\overline{p} + 8 = 10$, contradicting the hypothesis. So, $N(w_1) \cap V(I_b^*) = \emptyset$. Recalling also that $|I_i| = 3$ for each $i \in \{1, 2, ..., s\} \backslash \{a, b\}$, we conclude by Lemma 2(a1) that $N(w_1) \cap V(I_i^*) = \emptyset$ for each $i \in \{1, 2, ..., s\} \backslash \{a\}$, that is $N(w_1) \subseteq \{\xi_1, \xi_2, ..., \xi_s, w_2, w_3\} \backslash \{\xi_{a+1}, \xi_b\}$. Claim 2 is proved.  $\Delta$

If $\xi_{a+1} \neq \xi_b$ then by Claim 2, $|N(w_1)| \leq s = \delta - 1$, contradicting the fact that $|N(w_1)| \geq \delta$. Hence, $\xi_{a+1} = \xi_b$. Assume w.l.o.g. that $a = 1$.

**Claim 3.** $s = 2$.

**Proof.** Assume the contrary, that is $s \geq 3$. By Claim 1(2), $|I_3| = 3$. Put $I_3 = \xi_3 w_8 w_9 \xi_4$. If $\Upsilon(I_2, I_3) \neq \emptyset$ then we can argue as in Case 1.2.1.1. Otherwise, by Lemma 2(a1), $\Upsilon(I_3, I_i) = \emptyset$ for each $i \in \{1, 2, ..., s\} \backslash \{3\}$. If $w_9 \xi_2 \in E(G)$ then

$$\xi_1 \overrightarrow{C} w_3 z \overrightarrow{C} w_9 \xi_2 x_1 x_2 \xi_4 \overrightarrow{C} \xi_1$$

is longer than $C$, a contradiction. Hence, $N(w_9) \subseteq \{\xi_1, \xi_2, ..., \xi_s, w_8\} \backslash \{\xi_2\}$, contradicting the fact that $|N(w_9)| \geq \delta = s + 1$. Claim 3 is proved.  $\Delta$

By Claim 3,

$$\delta = s + 1 = 3, \ |C| = 9, \ n = 11, \ V(G) = V(C \cup P).$$

**Claim 4.** $N(w_1, w_2) = \{\xi_1, w_3\}$.

**Proof.** By Claims 2 and 3, $N(w_1) = \{\xi_1, w_2, w_3\}$, implying that $w_1 w_3 \in E(G)$. If $w_2 \xi_2 \in E(G)$ then

$$\xi_1 x_1 x_2 \xi_2 w_2 w_1 w_3 z \overrightarrow{C} \xi_1$$

is longer than $C$, a contradiction. Next, if $w_2 w_4 \in E(G)$ then

$$\xi_1 x_1 x_2 \xi_2 \overleftarrow{C} w_2 w_4 \overrightarrow{C} \xi_1$$

is longer than $C$, a contradiction. Further, if $w_2 w_5 \in E(G)$ then

$$\xi_1 x_1 x_2 \xi_2 w_3 w_1 w_2 w_5 \overrightarrow{C} \xi_1$$

is longer than $C$, a contradiction. Finally, if $w_2v \in E(G)$ for some $v \in \{w_6, w_7\}$ then $yz$ and $w_2v$ are independent intermediate edges between $I_1$ and $I_2$. By Lemma 2(a3), we reach a contradiction. Thus, $N(w_2) = \{\xi_1, w_1, w_3\}$ and Claim 4 follows.  $\Delta$

**Case 1.2.2.2.1.1.** $z = w_5$.
**Claim 5.** $N(w_6, w_7) = \{\xi_1, w_5\}$.
**Proof.** If $vw_3 \in E(G)$ for some $v \in \{w_6, w_7\}$ then

$$\xi_1 \overrightarrow{C} w_3 v \overleftarrow{C} \xi_2 x_1 x_2 \xi_1$$

is longer than $C$, a contradiction. Next, if $w_7\xi_2 \in E(G)$ then

$$\xi_1 \overrightarrow{C} w_3 w_5 w_6 w_7 \xi_2 x_2 x_1 \xi_1$$

is longer than $C$, a contradiction. Further, if $w_7w_4 \in E(G)$ then

$$\xi_1 \overrightarrow{C} w_3 w_5 w_6 w_7 w_4 \xi_2 x_2 x_1 \xi_1$$

is longer than $C$, a contradiction. Recalling also Claim 4, we get $N(w_7) = \{\xi_1, w_5, w_6\}$, implying that $w_7w_5 \in E(G)$. Now consider the neighborhood of $w_6$. If $w_6\xi_2 \in E(G)$ then

$$\xi_1 \overrightarrow{C} w_3 w_5 w_7 w_6 \xi_2 x_2 x_1 \xi_1$$

is longer than $C$, a contradiction. Next, if $w_6w_4 \in E(G)$ then

$$\xi_1 \overrightarrow{C} w_3 w_5 w_6 w_4 \xi_2 x_2 x_1 \xi_1$$

is longer than $C$, a contradiction. Recalling also Claim 4, we get $N(w_6) = \{\xi_1, w_5, w_7\}$ and Claim 5 follows.  $\Delta$

If $w_3w_4 \notin E(G)$ then by Claims 4 and 5, $G\backslash\{\xi_1, \xi_2, w_5\}$ has at least four connected components with vertex sets $\{x_1, x_2\}$, $\{w_1, w_2\}$, $\{w_6, w_7\}$ and $\{w_4\}$, contradicting the fact that $\tau \geq 1$. Hence, $w_3w_4 \in E(G)$ and it is easy to see that $G \in \Re_2$.

**Case 1.2.2.2.1.2.** $z = w_4$.
If $w_5w_3 \in E(G)$ then we can argue as in Case 1.2.2.2.1.1. Let $w_5w_3 \notin E(G)$.

**Case 1.2.2.2.1.2.1.** $w_5\xi_2 \notin E(G)$.
**Claim 6.** $N(w_5, w_6, w_7) = \{\xi_1, w_4\}$.
**Proof.** By the hypothesis, $N(w_5) \cap \{w_3, \xi_2\} = \emptyset$. Next, by Claim 4, $N(w_5) \cap \{w_1, w_2\} = \emptyset$. So, $N(w_5) \subseteq \{w_4, w_6, w_7, \xi_1\}$. Further, we can state (exactly as in Case 1.2.2.2.1.1) that $w_7w_3, w_7\xi_2, w_6w_3 \notin E(G)$. Moreover, if $w_6\xi_2 \in E(G)$ then

$$\xi_1 \overrightarrow{C} w_3 w_4 w_5 w_6 \xi_2 x_2 x_1 \xi_1$$

is longer than $C$, a contradiction. Recalling also Claim 4, we get $N(w_6, w_7) \subseteq \{w_4, w_5, \xi_1\}$ and Claim 6 follows. $\triangle$

If $N(w_6, w_7) = \{\xi_1, w_5\}$ then by Claims 4 and 6, $N(x_1, x_2) = \{\xi_1, \xi_2\}$, $N(w_1, w_2) = \{\xi_1, w_3\}$ and therefore, $G \in \Re_2$. Otherwise, $N(w_6, w_7) = \{\xi_1, w_4\}$, implying that $G \in \Re_1$.

**Case 1.2.2.2.1.2.2.** $w_5 \xi_2 \in E(G)$.
**Claim 7.** $N(w_6, w_7) = \{\xi_1, w_5\}$.
**Proof.** We can state (exactly as in Case 1.2.2.2.1.2.1) that $w_7 w_3, w_7 \xi_2, w_6 w_3 \notin E(G)$. If $w_7 w_4 \in E(G)$ then

$$\xi_1 \overrightarrow{C} w_3 w_4 w_7 w_6 w_5 \xi_2 x_2 x_1 \xi_1$$

is longer than $C$, a contradiction. So, using Claim 4, we conclude that $N(w_7) = \{w_5, w_6, \xi_1\}$. As for $N(w_6)$, if $w_6 \xi_2 \in E(G)$ then

$$\xi_1 \overrightarrow{C} w_3 w_4 w_5 w_6 \xi_2 x_2 x_1 \xi_1$$

is longer than $C$, a contradiction. If $w_6 w_4 \in E(G)$ then

$$\xi_1 \overrightarrow{C} w_3 w_4 w_6 w_5 \xi_2 x_2 x_1 \xi_1$$

is longer than $C$, a contradiction. Recalling also Claim 4, we can state that $N(w_6) = \{w_5, w_7, \xi_1\}$. Combining this with $N(w_7) = \{w_5, w_6, \xi_1\}$, we get $N(w_6, w_7) = \{w_5, \xi_1\}$. $\triangle$

By Claims 4 and 7,

$$N(x_1, x_2) = \{\xi_1, \xi_2\}, \ N(w_1, w_2) = \{\xi_1, w_3\}, \ N(w_6, w_7) = \{w_5, \xi_1\},$$

implying that $G \in \Re_2$.

**Case 1.2.2.2.2.** $y = w_2$.
If $z = w_4$ then $|C_1| > |C|$, a contradiction. If $z = w_7$ then $|C_2| > |C|$, again a contradiction. Hence, $z \in \{w_5, w_6\}$. Assume w.l.o.g. that $z = w_5$. Further, if $\Upsilon(I_1, I_2, ..., I_s) \neq \Upsilon(I_a, I_b)$ then by Lemma 2(a1), $\Upsilon(I_b, I_c) \neq \emptyset$ for some $c \in \{1, 2, ..., s\} \setminus \{a, b\}$. Since by Claim 2(2), $|I_c| = 3$, we can argue as in Case 1.2.1.1. Now let $\Upsilon(I_1, I_2, ..., I_s) = \Upsilon(I_a, I_b)$. If $w_1 w_3 \in E(G)$ then

$$\xi_a x_1 x_2 \xi_b \overleftarrow{C} w_3 w_1 w_2 w_5 \overrightarrow{C} \xi_a$$

is longer than $C$, a contradiction. Let $w_1 w_3 \notin EG)$. If $\Upsilon(I_a, I_b) = \{w_2 w_5\}$ then $G \setminus \{\xi_1, \xi_2, ..., \xi_s, w_2\}$ has at least $s + 2$ connected components, contradicting the fact that $\tau \geq 1$. Let $\Upsilon(I_a, I_b) = \{w_2 w_5, uv\}$. If $w_2 w_5$ and $uv$ are independent then by Lemma 2(a3), $|I_a| + |I_b| \geq 2\overline{p} + 8 = 10$, contradicting the hypothesis. Now let $w_2 w_5$ and $uv$ have a vertex in common. If $v = w_5$ then $u \in \{w_1, w_3\}$ and we can argue as in Case 1.2.2.2.1. Otherwise $u = w_2$ and hence $G\backslash\{\xi_1, \xi_2, ..., \xi_s, w_2\}$ has at least $s+2$ connected components, contradicting the fact that $\tau \geq 1$.

**Case 2.** $2 \leq \overline{p} \leq \delta - 3$.

It follows that $|N_C(x_i)| \geq \delta - \overline{p} \geq 3$ ($i = 1, 2$). If $N_C(x_1) \neq N_C(x_2)$ then by Lemma 1, $|C| \geq 4\delta - 2\overline{p} \geq 3\delta - \overline{p} + 3$, contradicting (1). Hence $N_C(x_1) = N_C(x_2)$, implying that $|I_i| \geq \overline{p} + 2$ ($i = 1, 2, ..., s$). Clearly, $s \geq |N_C(x_1)| \geq \delta - \overline{p} \geq 3$. If $s \geq \delta - \overline{p} + 1$ then

$$|C| \geq s(\overline{p} + 2) \geq (\delta - \overline{p} + 1)(\overline{p} + 2)$$

$$= (\delta - \overline{p} - 1)(\overline{p} - 1) + 3\delta - \overline{p} + 1 \geq 3\delta - \overline{p} + 3,$$

again contradicting (1). Hence $s = \delta - \overline{p} \geq 3$. It means that $x_1 x_2 \in E(G)$, that is $G[V(P)]$ is hamiltonian. By symmetric arguments, we can assume that $N_C(y) = N_C(x_1)$ for each $y \in V(P)$. If $\Upsilon(I_1, I_2, ..., I_s) = \emptyset$ then clearly $\tau < 1$, contradicting the hypothesis. Otherwise $\Upsilon(I_a, I_b) \neq \emptyset$ for some elementary segments $I_a$ and $I_b$. By definition, there is an intermediate path $L$ between $I_a$ and $I_b$. If $|L| \geq 2$ then by lemma 2(a1),

$$|I_a| + |I_b| \geq 2\overline{p} + 2|L| + 4 \geq 2\overline{p} + 8.$$

Hence

$$|C| = |I_a| + |I_b| + \sum_{i \in \{1,...,s\} \backslash \{a,b\}} |I_i| \geq 2\overline{p} + 8 + (s - 2)(\overline{p} + 2)$$

$$= (\delta - \overline{p} - 2)(\overline{p} - 1) + 3\delta - \overline{p} + 2 \geq 3\delta - \overline{p} + 3,$$

contradicting (1). Thus, $|L| = 1$, and therefore $\Upsilon(I_1, I_2, ..., I_s) \subseteq E(G)$. Let $L = yz$, where $y \in V(I_a^*)$ and $z \in V(I_b^*)$. Put

$$C_1 = \xi_a x_1 \overrightarrow{P} x_2 \xi_b \overleftarrow{C} yz \overrightarrow{C} \xi_a, \quad C_2 = \xi_a \overrightarrow{C} yz \overleftarrow{C} \xi_{a+1} x_1 \overrightarrow{P} x_2 \xi_{b+1} \overrightarrow{C} \xi_a.$$

By Lemma 2(a1),

$$|I_a| + |I_b| \geq 2\overline{p} + 2|L| + 4 = 2\overline{p} + 6,$$

which yields

$$|C| = |I_a| + |I_b| + \sum_{i \in \{1,...,s\} \backslash \{a,b\}} |I_i| \geq 2\overline{p} + 6 + (s - 2)(\overline{p} + 2)$$

$$= (s - 2)(\overline{p} - 2) + (\delta - \overline{p} - 4) + (3\delta - \overline{p} + 2).$$

Since $s = \delta - \overline{p} \geq 3$, we can state that if either $\overline{p} \geq 3$ or $\delta - \overline{p} \geq 4$ then $|C| \geq 3\delta - \overline{p} + 2$, contradicting (1). Otherwise $\overline{p} = 2$ and $\delta - \overline{p} = 3$, implying that $s = 3$ and $\delta = 5$. Assume w.l.o.g. that $a = 1$ and $b = 2$, i.e. $|I_1| + |I_2| \geq 10$, $|I_3| \geq 4$ and $|C| \geq 14$. On the other hand, by (1), $|C| \leq 3\delta - \overline{p} + 1 = 14$, which yields

$$n = 17, \ \delta = 5, \ |I_1| + |I_2| = 10, \ |I_3| = 4, \ |C| = 14, \ \overline{p} = 2. \tag{6}$$

If $|\Upsilon(I_1, I_2)| \geq 2$ then by Lemma 2(a2), $|I_1|+|I_2| \geq 2\overline{p}+7 = 11$, contradicting (6). Hence,
$$\Upsilon(I_1, I_2) = \{yz\}. \tag{7}$$

Since $|I_1| + |I_2| = 10$, we can assume w.l.o.g. that either $|I_1| = |I_2| = 5$ or $|I_1| = 4$, $|I_2| = 6$.

**Case 2.1.** $|I_1| = |I_2| = 5$.

Put $I_1 = \xi_1 w_1 w_2 w_3 w_4 \xi_2$ and $I_2 = \xi_2 w_5 w_6 w_7 w_8 \xi_3$. If $\Upsilon(I_i, I_3) \neq \emptyset$ for some $i \in \{1, 2\}$, then by Lemma 2(a1), $|I_i|+|I_3| \geq 2\overline{p}+6 = 10$. This implies $|I_3| \geq 5$, contradicting (6). Using also (7), we get
$$\Upsilon(I_1, I_2, I_3) = \{yz\}. \tag{8}$$

**Case 2.1.1.** $y \in \{w_2, w_3\}$.

Assume w.l.o.g. that $y = w_3$. If $z = w_5$ then $|C_1| > |C|$, a contradiction. Next, if $z \in \{w_7, w_8\}$ then $|C_2| > |C|$, again a contradiction. Hence $z = w_6$. Further, if $w_4 v \in E(G)$ for some $v \in \{w_1, w_2\}$ then

$$\xi_1 x_1 \overrightarrow{P} x_2 \xi_2 w_4 v \overrightarrow{C} w_3 w_6 \overrightarrow{C} \xi_1$$

is longer than $C$, a contradiction. Therefore, by (8), $N(w_4) \subseteq \{\xi_1, \xi_2, \xi_3, w_3\}$, contradicting the fact that $|N(w_4)| \geq \delta = 5$.

**Case 2.1.2.** $y \in \{w_1, w_4\}$.

Assume w.l.o.g. that $y = w_4$. Put $I_3 = \xi_3 w_9 w_{10} w_{11} \xi_1$. If $z \in \{w_6, w_7, w_8\}$ then $|C_2| > |C|$, a contradiction. Hence $z = w_5$. If $w_9 \xi_2 \in E(G)$ then

$$\xi_1 \overrightarrow{C} w_4 w_5 \overrightarrow{C} \xi_3 x_2 \overleftarrow{P} x_1 \xi_2 w_9 \overrightarrow{C} \xi_1$$

is longer than $C$, a contradiction. Hence, $w_9 \xi_2 \notin E(G)$. By (8), $N(w_9) \subseteq \{\xi_1, \xi_3, w_{10}, w_{11}\}$, contradicting the fact that $|N(w_9)| \geq \delta = 5$.

**Case 2.2.** $|I_1| = 4$ and $|I_2| = 6$.

Put $I_1 = \xi_1 w_1 w_2 w_3 \xi_2$ and $I_2 = \xi_2 w_4 w_5 w_6 w_7 w_8 \xi_3$. By Lemma 2(a1),
$$\Upsilon(I_1, I_3) = \emptyset. \tag{9}$$

**Case 2.2.1.** $y = w_2$.

If $z \in \{w_4, w_5\}$ then $|C_1| > |C|$, a contradiction. If $z \in \{w_7, w_8\}$ then $|C_2| > |C|$, again a contradiction. Hence, $z = w_6$. Next, if $w_1 w_3 \in E(G)$ then

$$\xi_1 x_1 \overrightarrow{P} x_2 \xi_2 w_3 w_1 w_2 w_6 \overrightarrow{C} \xi_1$$

is longer than $C$, a contradiction. By (7) and (9), $N(w_1) \subseteq \{\xi_1, \xi_2, \xi_3, w_2\}$, contradicting the fact that $|N(w_1)| \geq \delta = 5$.

**Case 2.2.2.** $y = w_1$.

If $z \in \{w_4, w_5, w_6\}$ then $|C_1| > |C|$, a contradiction. Next, if $z = w_8$ then $|C_2| > |C|$, a contradiction. Hence, $z = w_7$. Further, if $w_3\xi_1 \in E(G)$ then

$$\xi_1 w_3 w_2 w_1 w_7 \overleftarrow{C} \xi_2 x_1 \overrightarrow{P} x_2 \xi_3 \overrightarrow{C} \xi_1$$

is longer than $C$, a contradiction. By (7) and (9), $N(w_3) \subseteq \{\xi_2, \xi_3, w_1, w_2\}$, contradicting the fact that $|N(w_3)| \geq \delta = 5$.

**Case 2.2.3.** $y = w_3$.

If $z = w_4$ then $|C_1| > |C|$, a contradiction. Next, if $z \in \{w_6, w_7, w_8\}$ then $|C_2| > |C|$, a contradiction. Hence, $z = w_5$. Further, if $w_1\xi_2 \in E(G)$ then

$$\xi_1 x_1 \overrightarrow{P} x_2 \xi_2 w_1 w_2 w_3 w_5 \overrightarrow{C} \xi_1$$

is longer than $C$, a contradiction. By (7) and (9), $N(w_1) \subseteq \{\xi_1, \xi_3, w_2, w_3\}$, contradicting the fact that $|N(w_1)| \geq \delta = 5$.

**Case 3.** $2 \leq \overline{p} = \delta - 2$.

It follows that $|N_C(x_i)| \geq \delta - \overline{p} = 2$ $(i = 1, 2)$ and

$$\delta \geq 4. \tag{10}$$

If $N_C(x_1) \neq N_C(x_2)$ then by Lemma 1, $|C| \geq 4\delta - 2\overline{p} = 3\delta - \overline{p} + 2$, contradicting (1). Hence, $N_C(x_1) = N_C(x_2)$. Clearly, $s = |N_C(x_1)| \geq 2$. Further, if $s \geq 3$ then

$$|C| \geq s(\overline{p} + 2) \geq 3\delta \geq 3\delta - \overline{p} + 2,$$

again contradicting (1). So, $s = 2$. It follows that $x_1 x_2 \in E(G)$, that is $G[V(P)]$ is hamiltonian. By symmetric arguments, $N_C(v) = N_C(x_1) = \{\xi_1, \xi_2\}$ for each $v \in V(P)$. If $\Upsilon(I_1, I_2) = \emptyset$ then clearly $\tau < 1$, contradicting the hypothesis. Otherwise, there is an intermediate path $L = y\overrightarrow{L}z$ such that $y \in V(I_1^*)$ and $z \in V(I_2^*)$. If $|L| \geq 2$ then by Lemma 2(a1),

$$|C| = |I_1| + |I_2| \geq 2\overline{p} + 2|L| + 4 \geq 2\overline{p} + 8 = 3\delta - \overline{p} + 2,$$

contradicting (1). Hence $|L| = 1$, that is $L = yz$, implying that $\Upsilon(I_1, I_2) \subseteq E(G)$. If $|\Upsilon(I_1, I_2)| \geq 3$ then by Lemma 2(a2), $|C| = |I_1| + |I_2| \geq 2\overline{p} + 8 = 3\delta - \overline{p} + 2$, contradicting (1). Hence $1 \leq |\Upsilon(I_1, I_2)| \leq 2$. Put

$$Q = \xi_1 x_1 \overrightarrow{P} x_2 \xi_2.$$

**Case 3.1.** $|\Upsilon(I_1, I_2)| = 1$.

It follows that

$$\Upsilon(I_1, I_2) = \{yz\}. \tag{11}$$

By Lemma 2(a1), $|C| = |I_1| + |I_2| \geq 2\overline{p} + 2|L| + 4 = 2\delta + 2$. On the other hand, by (1), $|C| \leq 2\delta + 3$. So,

$$2\delta + 2 \leq |C| \leq 2\delta + 3. \tag{12}$$

Put
$$C_1 = \xi_1 \overrightarrow{Q} \xi_2 \overrightarrow{C} zy \overleftarrow{C} \xi_1, \quad C_2 = \xi_1 \overrightarrow{Q} \xi_2 \overleftarrow{C} yz \overrightarrow{C} \xi_1,$$
$$|\xi_1 \overrightarrow{C} y| = d_1, \quad |y \overrightarrow{C} \xi_2| = d_2, \quad |\xi_2 \overrightarrow{C} z| = d_3, \quad |z \overrightarrow{C} \xi_1| = d_4.$$

**Case 3.1.1.** $d_i = 1$ for some $i \in \{1, 2, 3, 4\}$.

Assume w.l.o.g. that $d_1 = 1$. Since $|I_i| \geq \overline{p} + 2 = \delta$ $(i = 1, 2)$, we have $d_2 \geq \delta - 1$.

**Case 3.1.1.1.** $d_2 = \delta - 1$.

It follows that $|I_1| = \delta$. By (12), $|I_2| \geq \delta + 2$. Next, by (11), $N(\xi_2^-) \subset V(I_1)$. Recalling that $|I_1| = \delta$, we have $N(\xi_2^-) = V(I_1) \backslash \{\xi_2^-\}$. In particular, $\xi_2^- \xi_1 \in E(G)$. If $d_3 \geq 3$ then
$$|C| \geq |\xi_1 \overrightarrow{Q} \xi_2 \overrightarrow{C} zy \overrightarrow{C} \xi_2^- \xi_1| \geq (\overline{p} + 2) + d_3 + d_2 + 2 \geq 2\delta + 4,$$
contradicting (12). Let $d_3 \leq 2$. Further, since $|I_2| \geq \delta + 2$, we have $d_4 \geq \delta$. Hence,
$$|C| \geq |C_2| \geq (\overline{p} + 2) + d_2 + d_4 + 1 \geq 3\delta.$$
By (12), $3\delta \leq |C| \leq 2\delta + 3$, which yields $\delta \leq 3$, contradicting (10).

**Case 3.1.1.2.** $d_2 = \delta$.

If $d_4 \geq 3$ then
$$|C| \geq |C_2| \geq \delta + d_2 + d_4 + 1 \geq 2\delta + 4,$$
contradicting (12). Let $d_4 \leq 2$. Since $|I_1| = d_2 + 1 = \delta + 1$, we have by (12), $|I_2| \geq \delta + 1$, that is $d_3 \geq \delta - 1$.

**Case 3.1.1.2.1.** $d_3 = \delta - 1$.

If $d_4 = 1$ then we can reach a contradiction as in Case 3.1.1.1. Let $d_4 = 2$. If $\xi_2^+ \xi_1 \in E(G)$ then
$$|C| \geq |\xi_1 \overrightarrow{Q} \xi_2 \overleftarrow{C} yz \overleftarrow{C} \xi_2^+ \xi_1| \geq 3\delta + 1.$$
By (12), $3\delta + 1 \leq |C| \leq 2\delta + 3$, that is $\delta \leq 2$, contradicting (10). Further, if $\xi_2^+ z^+ \in E(G)$ then
$$|C| \geq |\xi_1 z^+ \xi_2^+ \overrightarrow{C} zy \overrightarrow{C} \xi_2 \overleftarrow{Q} \xi_1| \geq 3\delta + 2.$$
By (12), $3\delta + 2 \leq |C| \leq 2\delta + 3$, which yields $\delta \leq 1$, contradicting (10). Therefore, by (11), $N(\xi_2^+) \subseteq V(\xi_2 \overrightarrow{C} z) \backslash \{\xi_2^+\}$, implying that $|N(\xi_2^+)| \leq \delta - 1$. This contradicts the fact that $d(\xi_2^+) \geq \delta$.

**Case 3.1.1.2.2.** $d_3 = \delta$.
**Case 3.1.1.2.2.1.** $d_4 = 2$.
If $z^+ w \in E(G)$ for some $w \in V(\xi_2^+ \overrightarrow{C} z^-)$, then
$$|C| \geq |\xi_1 z^+ w \overrightarrow{C} zy \overrightarrow{C} \xi_2 \overleftarrow{Q} \xi_1| \geq 2\delta + 4,$$

contradicting (12). Therefore, by (11), $N(z^+) \subseteq \{\xi_1, \xi_2, z\}$, contradicting (10).

**Case 3.1.1.2.2.2**. $d_4 = 1$.
Assume that $\xi_1 w \in E(G)$ for some $w \in V(y^+ \overrightarrow{C} \xi_2^-)$. If $w \neq y^+$ then

$$|C| \geq |\xi_1 w \overleftarrow{C} yz \overleftarrow{C} \xi_2 \overleftarrow{Q} \xi_1| \geq 2\delta + 4,$$

contradicting (12). Let $w = y^+$. If $\xi_2^- y \in E(G)$ then

$$|C| \geq |\xi_1 y^+ \overrightarrow{C} \xi_2^- yz \overleftarrow{C} \xi_2 \overleftarrow{Q} \xi_1| = 3\delta + 1 > 2\delta + 4,$$

contradicting (12). Let $\xi_2^- y \notin E(G)$. Further, if $\xi_2^- \xi_1 \in E(G)$ then

$$|C| \geq |\xi_1 \xi_2^- \overleftarrow{C} yz \overleftarrow{C} \xi_2 \overleftarrow{Q} \xi_1| \geq 3\delta + 1 > 2\delta + 4,$$

contradicting (12). Therefore, by (11), $N(\xi_2^-) \subseteq V(y^+ \overrightarrow{C} \xi_2) \backslash \{\xi_2^-\}$, implying that $|N(\xi_2^-)| \leq \delta - 1$, a contradiction. So, $\xi_1 w \notin E(G)$ for each $w \in V(y^+ \overrightarrow{C} \xi_2^-)$. By a symmetric argument, $\xi_1 w \notin E(G)$ for each $w \in V(\xi_2^+ \overrightarrow{C} z^-)$. So, $G \in \Re_1$.

**Case 3.1.1.2.3**. $d_3 \geq \delta + 1$.
If either $d_3 \geq \delta + 2$ or $d_4 \geq 2$ then $|I_2| \geq \delta + 3$ and hence $|C| = |I_1| + |I_2| \geq 2\delta + 4$, contradicting (12). Hence, $d_3 = \delta + 1$, $d_4 = 1$. If $N(z) \cap V(\xi_2^+ \overrightarrow{C} z^{-2}) = \emptyset$ then $\delta = 4$ and $G \in \Re_3$. Otherwise, $G \in \Re_1$.

**Case 3.1.1.3**. $d_2 \geq \delta + 1$.
If either $d_2 \geq \delta + 2$ or $d_4 \geq 2$ then $|C| \geq |C_2| \geq 2\delta + 4$, contradicting (12). Otherwise $d_2 = \delta + 1$ and $d_4 = 1$. By (12), $d_3 \leq \delta$. If $d_3 = \delta$ then we can argue as in Case 3.1.1.2.3. Let $d_3 = \delta - 1$. If $\xi_2^+ \xi_1 \in E(G)$ then

$$|C| \geq |\xi_1 \xi_2^+ \overrightarrow{C} zy \overrightarrow{C} \xi_2 \overleftarrow{Q} \xi_1| \geq 3\delta + 1 > 2\delta + 4,$$

a contradiction. Therefore, by (11), $N(\xi_2^+) \subseteq V(\xi_2 \overrightarrow{C} z) \backslash \{\xi_2^+\}$, implying that $|N(\xi_2^+)| \leq \delta - 1$, a contradiction.

**Case 3.1.2**. $d_i \geq 2$ ($i = 1, 2, 3, 4$).
If $d_i \geq \delta - 1$ ($i = 1, 2, 3, 4$) then $|C| \geq 4\delta - 4 \geq 2\delta + 4$, contradicting (12). Assume w.l.o.g. that $d_1 \leq \delta - 2$. If $N(y^-) \subseteq V(\xi_1 \overrightarrow{C} y) \cup \{\xi_2\}$ then clearly $|N(y^-)| \leq \delta - 1$, a contradiction. Otherwise, by (11), $y^- w \in E(G)$ for some $w \in V(y^+ \overrightarrow{C} \xi_2^-)$.

**Case 3.1.2.1**. $w \neq y^+$.
Clearly,
$$|C| \geq |C_2| \geq \delta + d_2 + d_4 + 1,$$
$$|C| \geq |\xi_1 \overrightarrow{C} y^- w \overleftarrow{C} yz \overleftarrow{C} \xi_2 \overleftarrow{Q} \xi_1| \geq \delta + d_1 + d_3 + 3.$$

By summing, we get $2|C| \geq 2\delta + \sum_{i=1}^{4} d_i + 4 = 2\delta + |C| + 4$, that is $|C| \geq 2\delta + 4$, contradicting (12).

**Case 3.1.2.2.** $w = y^+$.

We have
$$|C| \geq |\xi_1\overrightarrow{Q}\xi_2\overleftarrow{C}y^+y^-yz\overrightarrow{C}\xi_1| \geq \delta + d_2 + d_4 + 2,$$
$$|C| \geq |\xi_1\overrightarrow{C}y^-y^+yz\overleftarrow{C}\xi_2\overleftarrow{Q}\xi_1| \geq \delta + d_1 + d_3 + 2.$$

By summing, we get $2|C| \geq 2\delta + \sum_{i=1}^{4} d_i + 4 = 2\delta + |C| + 4$, that is $|C| \geq 2\delta + 4$, contradicting (12).

**Case 3.2.** $|\Upsilon(I_1, I_2)| = 2$.

If $\Upsilon(I_1, I_2)$ consists of two independent edges, then by Lemma 2(a3), $|C| = |I_1| + |I_2| \geq 2\overline{p} + 8 = 2\delta + 4$, contradicting (12). Otherwise we can assume w.l.o.g. that
$$\Upsilon(I_1, I_2) = \{yz_1, yz_2\}, \tag{13}$$
where $y \in V(I_1^*)$ and $z_1, z_2 \in V(I_2^*)$. Assume w.l.o.g. that $z_1, z_2$ occur on $\xi_2\overrightarrow{C}\xi_1$ in this order. Put
$$C_1 = \xi_1\overrightarrow{Q}\xi_2\overrightarrow{C}z_2y\overleftarrow{C}\xi_1, \quad C_2 = \xi_1\overrightarrow{Q}\xi_2\overrightarrow{C}yz_1\overrightarrow{C}\xi_1,$$
$$|\xi_1\overrightarrow{C}y| = d_1, \ |y\overrightarrow{C}\xi_2| = d_2, \ |\xi_2\overrightarrow{C}z_1| = d_3, \ |z_2\overrightarrow{C}\xi_1| = d_4, \ |z_1\overrightarrow{C}z_2| = d_5.$$

By Lemma 2(a2), $|C| = |I_1| + |I_2| \geq 2\overline{p} + 7 = 2\delta + 3$. On the other hand, by (1), $|C| \leq 3\delta - \overline{p} + 1 = 2\delta + 3$, implying that
$$|C| = 2\delta + 3. \tag{14}$$

Assume first that $d_5 \geq 2$. Clearly,
$$|C| \geq |C_2| \geq \delta + d_2 + d_4 + 3,$$
$$|C| \geq |C_1| = \delta + d_1 + d_3 + d_5 + 1.$$

By summing, we get $2|C| \geq 2\delta + \sum_{i=1}^{5} d_i + 4 = 2\delta + |C| + 4$, that is $|C| \geq 2\delta + 4$, contradicting (14). Now let $d_5 = 1$.

**Case 3.2.1.** $d_i = 1$ for some $i \in \{1, 2, 3, 4\}$.

We can assume w.l.o.g. that either $d_1 = 1$ or $d_4 = 1$.

**Case 3.2.1.1.** $d_1 = 1$.

Since $|I_i| \geq \delta$ ($i = 1, 2$), we have $d_2 \geq \delta - 1$.

**Case 3.2.1.1.1.** $d_2 = \delta - 1$.

It follows that $|I_1| = \delta$. By (14), $|I_2| = \delta + 3$. Further, by (13), $N(\xi_2^-) \subset V(I_1)$. Since $|I_1| = \delta$, we have $N(\xi_2^-) = V(I_1) \setminus \{\xi_2^-\}$, implying that $\xi_2^- \xi_1 \in E(G)$. If $d_3 \geq 3$ then
$$|C| \geq |\xi_1\xi_2^-\overleftarrow{C}yz_2\overleftarrow{C}\xi_2\overleftarrow{Q}\xi_1| \geq 2\delta + 4,$$

contradicting (14). Hence, $d_3 \leq 2$, implying that $d_4 = |I_2| - d_3 - d_5 \geq \delta$. But then
$$|C| \geq |C_2| = 3\delta + 1 > 2\delta + 4,$$
contradicting (14).

**Case 3.2.1.1.2**. $d_2 = \delta$.
It follows that $|I_1| = \delta + 1$. By (14), $|I_2| = \delta + 2$. If $d_4 \geq 2$ then
$$|C| \geq |C_2| \geq 2\delta + 4,$$
contradicting (14). Hence, $d_4 = 1$ and $d_3 = \delta$. If $z_2 w \in E(G)$ for some $w \in V(\xi_2^+ \overrightarrow{C} z_1^-)$ then
$$|C| \geq |\xi_1 \overleftarrow{C} z_2 w \overrightarrow{C} z_1 y \overrightarrow{C} \xi_2 \overleftarrow{Q} \xi_1| \geq 2\delta + 4,$$
contradicting (14). Therefore, by (13), $N(z_2) \subseteq \{\xi_1, \xi_2, y, z_1\}$, implying that $\delta = 4$. Next, if $\xi_1 w \in E(G)$ for some $w \in V(y^+ \overrightarrow{C} \xi_2^-)$ then
$$|C| \geq |\xi_1 w \overleftarrow{C} y z_2 \overleftarrow{C} \xi_2 \overleftarrow{Q} \xi_1| \geq 2\delta + 4,$$
contradicting (14). Further, if $\xi_1 w \in E(G)$ for some $w \in V(\xi_2^+ \overrightarrow{C} z_1^-)$ then
$$|C| \geq |\xi_1 w \overrightarrow{C} z_2 y \overrightarrow{C} \xi_2 \overleftarrow{Q} \xi_1| \geq 2\delta + 4,$$
contradicting (14). Thus, by (13), $N(\xi_1) \subseteq V(P) \cup \{y, z_1, z_2, \xi_2\}$, implying that $G \in \Re_3$.

**Case 3.2.1.1.3**. $d_2 \geq \delta + 1$.
It follows that $|C| \geq |C_2| \geq 2\delta + 4$, contradicting (14).

**Case 3.2.1.2**. $d_4 = 1$.
Since $|I_2| \geq \delta$ and $d_5 = 1$, we have $d_3 \geq \delta - 2$.

**Case 3.2.1.2.1**. $d_3 = \delta - 2$.
It follows that $|I_2| = \delta$. By (14), $|I_1| = \delta + 3$. If $d_1 \geq 4$ then
$$|C| \geq |C_1| \geq 2\delta + 4,$$
contradicting (14). Let $d_1 \leq 3$ which yields $d_2 \geq \delta$. If $z_2 w \in E(G)$ for some $w \in V(\xi_2^+ \overrightarrow{C} z_1^-)$ then
$$|C| \geq |\xi_1 z_2 w \overrightarrow{C} z_1 y \overrightarrow{C} \xi_2 \overleftarrow{Q} \xi_1| \geq 2\delta + 4,$$
contradicting (14). Therefore, by (13), $N(z_2) \subseteq \{\xi_1, \xi_2, y, z_1\}$, implying that $\delta = 4$. But then $d_3 = 2$ and $N(\xi_2^+) \subseteq \{\xi_1, \xi_2, z_1\}$, contradicting the fact that $d(\xi_2^+) \geq \delta \geq 4$.

**Case 3.2.1.2.2**. $d_3 = \delta - 1$.

It follows that $|I_2| = \delta + 1$. By (13), $|I_1| = \delta + 2$. If $d_1 \geq 3$ then

$$|C| \geq |C_1| \geq 2\delta + 4,$$

contradicting (14). Let $d_1 \leq 2$, implying that $d_2 \geq \delta$. If $z_2 w \in E(G)$ for some $w \in V(\xi_2^+ \overrightarrow{C} z_1^-)$ then

$$|C| \geq |\xi_1 z_2 w \overrightarrow{C} z_1 y \overrightarrow{C} \xi_2 \overleftarrow{Q} \xi_1| \geq 2\delta + 4,$$

contradicting (14). Otherwise, by (13), $N(z_2) \subseteq \{\xi_1, \xi_2, y, z_1\}$, implying that $\delta = 4$ and $N(\xi_2^+) \subseteq \{\xi_1, \xi_2, z_1, z_1^-\}$. By (10), $\xi_2^+ \xi_1 \in E(G)$ and therefore,

$$|C| \geq |\xi_1 \xi_2^+ \overrightarrow{C} z_2 y \overrightarrow{C} \xi_2 \overleftarrow{Q} \xi_1| \geq 2\delta + 5,$$

contradicting (14).

**Case 3.2.1.2.3**. $d_3 = \delta$.

It follows that $|I_2| = \delta + 2$. By (14), $|I_1| = \delta + 1$. If $d_1 \geq 2$ then

$$|C| \geq |C_1| \geq 2\delta + 4,$$

contradicting (14). Hence $d_1 = 1$ and $d_2 = \delta$. If $z_2 w \in E(G)$ for some $w \in V(\xi_2^+ \overrightarrow{C} z_1^-)$ then

$$|C| \geq |\xi_1 z_2 w \overrightarrow{C} z_1 y \overrightarrow{C} \xi_2 \overleftarrow{Q} \xi_1| \geq 2\delta + 4,$$

contradicting (14). Otherwise, by (13), $N(z_2) \subseteq \{\xi_1, \xi_2, y, z_1\}$, implying that $\delta = 4$. Further, if $\xi_1 w \in E(G)$ for some $w \in V(y^+ \overrightarrow{C} \xi_2^-)$ then

$$|C| \geq |\xi_1 w \overleftarrow{C} y z_2 \overleftarrow{C} \xi_2 \overleftarrow{Q} \xi_1| \geq 2\delta + 4,$$

contradicting (14). Let $\xi_1 w \notin E(G)$ for each $w \in V(y^+ \overrightarrow{C} \xi_2^-)$. Finally, if $\xi_1 w \in E(G)$ for some $w \in V(\xi_2^+ \overrightarrow{C} z_1^-)$ then

$$|C| \geq |\xi_1 w \overrightarrow{C} z_2 y \overrightarrow{C} \xi_2 \overleftarrow{Q} \xi_1| \geq 2\delta + 4,$$

contradicting (14). Hence, $\xi_1 w \notin E(G)$ for each $w \in V(\xi_2^+ \overrightarrow{C} z_1^-)$, implying that $G \in \Re_3$.

**Case 3.2.1.2.4**. $d_3 \geq \delta + 1$.

By (14), $d_3 = \delta + 1$, $|I_2| = \delta + 3$ and $|I_1| = \delta$. This implies

$$|C| \geq |C_1| \geq 2\delta + 4,$$

contradicting (14).

**Case 3.2.2**. $d_i \geq 2$ $(i = 1, 2, 3, 4)$.

If $d_i \geq \delta - 1$ ($i = 1, 2, 3, 4$) then $|C| \geq 4\delta - 4 \geq 2\delta + 4$, contradicting (14). Otherwise we can assume w.l.o.g. that either $d_1 \leq \delta - 2$ or $d_4 \leq \delta - 2$.

**Case 3.2.2.1**. $d_1 \leq \delta - 2$.
If $N(y^-) \subseteq V(\xi_1 \overrightarrow{C} y) \cup \{\xi_2\}$ then $|N(y^-)| \leq \delta - 1$, a contradiction. Otherwise, by (13), $y^- w \in E(G)$ for some $w \in V(y^+ \overrightarrow{C} \xi_2^-)$. Clearly,

$$|C| \geq |C_2| \geq \delta + d_2 + d_4 + 2,$$

$$|C| \geq |\xi_1 \overrightarrow{C} y^- w \overleftarrow{C} y z_2 \overleftarrow{C} \xi_2 \overleftarrow{Q} \xi_1| \geq \delta + d_1 + d_3 + 3.$$

By summing, we get $2|C| \geq 2\delta + \sum_{i=1}^{5} d_i + 4 \geq 2\delta + |C| + 4$, that is $|C| \geq 2\delta + 4$, contradicting (14).

**Case 3.2.2.2**. $d_4 \leq \delta - 2$.
If either $d_1 \leq \delta - 2$ or $d_2 \leq \delta - 2$ then we can argue as in Case 3.2.2.1. Let $d_i \geq \delta - 1$ ($i = 1, 2$). If $d_3 \geq 3$ then

$$|C| \geq |C_1| \geq \delta + d_1 + d_3 + 2 \geq 2\delta + 4,$$

contradicting (14). Hence, $d_3 = 2$. By a symmetric argument, $d_4 = 2$. Next, if $z_2^+ z_1 \in E(G)$ then

$$|C| \geq |\xi_1 \overrightarrow{C} y z_2 z_2^+ z_1 z_1^- \xi_2 \overleftarrow{Q} \xi_1| \geq 2\delta + 4,$$

contradicting (14). Let $z_2^+ z_1 \notin E(G)$. Finally, if $z_2^+ z_1^- \in E(G)$ then

$$|C| \geq |\xi_1 \overrightarrow{C} y z_1 z_2 z_2^+ z_1^- \xi_2 \overleftarrow{Q} \xi_1| \geq 2\delta + 4,$$

contradicting (14). So, by (13), $N(z_2^+) \subseteq \{\xi_1, \xi_2, z_2\}$, implying that $d(z_2^+) \leq 3$, which contradicts (10).

**Case 4**. $2 \leq \overline{p} = \delta - 1$.
It follows that $s \geq |N_C(x_i)| \geq \delta - \overline{p} = 1$ ($i = 1, 2$). By (1),

$$|C| = 2\delta + 2. \tag{15}$$

**Case 4.1**. $s \geq 3$.
Clearly, there are at least two segments among $I_1, I_2, ..., I_s$ of length at least $\overline{p} + 2$. This yields $|C| > 2(\overline{p} + 2) = 2\delta + 2$, contradicting (15).

**Case 4.2**. $s = 2$.
Clearly, $|I_i| \geq \overline{p} + 2 = \delta + 1$. By (15), $|I_1| = |I_2| = \delta + 1$ and $V(G) = V(C \cup P)$. If $\Upsilon(I_1, I_2) \neq \emptyset$ then by Lemma 2(a1), $|I_1| + |I_2| \geq 2\overline{p} + 6 = 2\delta + 4$, contradicting (15). Let $\Upsilon(I_1, I_2) = \emptyset$. Further, if there is an edge $yz$ such that $y \in V(P)$ and $z \in V(I_1^*) \cup V(I_2^*)$, say $z \in V(I_1^*)$, then

$$|I_1| = |\xi_1 \overrightarrow{C} z| + |z \overrightarrow{C} \xi_2| \geq |\xi_1 x_1 \overrightarrow{P} y z| + |z y \overrightarrow{P} x_2 \xi_2| \geq |P| + 4 = \delta + 3,$$

a contradiction. Otherwise, $G\backslash\{\xi_1,\xi_2\}$ has three connected components, contradicting the fact that $\tau \geq 1$.

**Case 4.3**. $s = 1$.

It follows that $N_C(x_1) = N_C(x_2) = \{y_1\}$ for some $y_1 \in V(C)$. Since $G$ is 2-connected, there is an edge $zw$ such that $z \in V(P)$ and $w \in V(C)\backslash\{y_1\}$. Clearly, $z \not\in \{x_1, x_2\}$ and $x_2z^- \in E(G)$. If $z \in \{x_1, x_2\}$ then $s \geq 2$, a contradiction. Let $z \not\in \{x_1, x_2\}$. Further, since $\overline{p} = \delta - 1$ and $|N_C(x_2)| = 1$, we have $x_2z^- \in E(G)$. Then replacing $P$ with $x_1\overrightarrow{P}z^-x_2\overleftarrow{P}z$, we can argue as in Case 4.2. ∎

**Proof of Theorem 2**. The proof follows from Theorem 1 immediately, since $\kappa(G) = 2$ and $\tau(G) = 1$ for each $G \in \Re$. ∎

Institute for Informatics and Automation Problems
National Academy of Sciences
P. Sevak 1, Yerevan 0014, Armenia
E-mail: zhora@ipia.sci.am

# References

[1] D. Bauer and E. Schmeichel, Long cycles in tough graphs, Technical Report 8612, Stevens Institute of Technology, 1986.

[2] A. Bigalke and H.A. Jung, Über Hamiltonische Kreise und unabhängige Ecken in Graphen, Monatsh. Math. 88 (1979) 195-210.

[3] J.A. Bondy and U.S.R. Murty, Graph Theory with Applications. Macmillan, London and Elsevier, New York (1976).

[4] C.St.J.A. Nash-Williams, Edge-disjoint hamiltonian cycles in graphs with vertices of large valency, in: L. Mirsky, ed., "Studies in Pure Mathematics", pp. 157-183, Academic Press, San Diego/London (1971).

[5] Zh.G. Nikoghosyan, A Degree Condition for Dominating Cycles in t-tough graphs with $t > 1$, arXiv.1201.1551.